\documentclass{ieeeaccess}
\usepackage[utf8]{inputenc}
\usepackage{amsmath,amssymb,amsfonts}
\usepackage[sort]{natbib}
\usepackage{graphicx}
\usepackage{multirow}
\usepackage{textcomp}
\usepackage{bm}
\makeatletter
\AtBeginDocument{\DeclareMathVersion{bold}
\SetSymbolFont{operators}{bold}{T1}{times}{b}{n}
\SetSymbolFont{NewLetters}{bold}{T1}{times}{b}{it}
\SetMathAlphabet{\mathrm}{bold}{T1}{times}{b}{n}
\SetMathAlphabet{\mathit}{bold}{T1}{times}{b}{it}
\SetMathAlphabet{\mathbf}{bold}{T1}{times}{b}{n}
\SetMathAlphabet{\mathtt}{bold}{OT1}{pcr}{b}{n}
\SetSymbolFont{symbols}{bold}{OMS}{cmsy}{b}{n}
\renewcommand\boldmath{\@nomath\boldmath\mathversion{bold}}}
\makeatother

\def\BibTeX{{\rm B\kern-.05em{\sc i\kern-.025em b}\kern-.08em
    T\kern-.1667em\lower.7ex\hbox{E}\kern-.125emX}}

\newcommand{\mat}[1]{#1}
\renewcommand{\vec}[1]{#1}
\renewcommand{\thetable}{\Roman{table}}
\newcommand{\R}[1]{\mathbb{R}^{#1}}
\begin{document}
\history{Date of publication xxxx 00, 0000, date of current version xxxx 00, 0000.}
\doi{10.1109/ACCESS.2023.1120000}

\title{Learning Nonlinear Dynamics Using Kalman Smoothing}
\author{
    \uppercase{Jacob M. Stevens-Haas}\authorrefmark{1},
    \uppercase{Yash Bhangale}\authorrefmark{2},
    J. Nathan Kutz\authorrefmark{1}\authorrefmark{3},
    and Aleksandr Aravkin\authorrefmark{1}}
\address[1]{
    Department of Applied Mathematics, University of Washington, Seattle, WA 98195-3925 USA}
\address[2]{Department of Mechanical Engineering, University of Washington}
\address[3]{Department of Electrical and Computer Engineering, University of Washington}
\tfootnote{This work was supported in part by the US National Science Foundation (NSF) AI Institute for Dynamical Systems (dynamicsai.org), grant 2112085 and by the US Department of Veterans Affairs under the Post-9/11 GI Bill}

\markboth
{Stevens-Haas \headeretal: Learning Nonlinear Dynamics Using Kalman Smoothing}
{Stevens-Haas \headeretal: Learning Nonlinear Dynamics Using Kalman Smoothing}

\corresp{Corresponding author: Jacob M. Stevens-Haas (e-mail: jacob.stevens.haas@gmail.com).}

\begin{abstract}
    Identifying Ordinary Differential Equations (ODEs) from measurement data requires both fitting the dynamics and assimilating, either implicitly or explicitly, the measurement data.
    The {\em Sparse Identification of Nonlinear Dynamics} (SINDy) method involves a derivative estimation (and optionally, smoothing) step and a sparse regression step on a library of candidate ODE terms.
    Kalman smoothing is a classical framework for assimilating the measurement data with known noise statistics.
    Previously, derivatives in SINDy and its python package, pysindy, had been estimated by finite difference, L1 total variation minimization, or local filters like Savitzky-Golay.
    In contrast, Kalman allows discovering ODEs that best recreate the essential dynamics in simulation, even in cases when it does not perform as well at recovering coefficients, as measured by their F1 score and mean absolute error.
    We have incorporated Kalman smoothing, along with hyperparameter optimization, into the existing pysindy architecture, allowing for rapid adoption of the method.
    Numerical experiments on a number of dynamical systems show 
    Kalman smoothing to be the most amenable to parameter selection and best at preserving problem structure in the presence of noise.
\end{abstract}


\titlepgskip=-21pt

\maketitle
\frenchspacing
\section{Introduction}
The method of {\em Sparse Identification of Nonlinear Dynamics} (SINDy)~\citep{Brunton2016,brunton2022data} seeks to discover a differential or partial differential equation governing an arbitrary, temporally measured system.
The method takes as input some coordinate measurements over time, such as angles between molecular bonds~\citep{Boninsegna2018} or a spatial field, such as wave heights~\citep{Rudy2017}, and returns the best ordinary or partial differential equation (ODE or PDE) from a library of candidate terms.
However, the method struggles to accommodate significant measurement noise, which is typical of real-world systems.
On the other hand, Kalman theory~\citep{kalman,KalBuc} has a half-century history of assimilating measurement noise to smooth a trajectory, with well-studied and rigorously characterized noise properties~\citep{welch1995introduction}.
We integrate the mature and well-established theory of Kalman with the emerging SINDy technology and combine with generalized cross validation (GCV) parameter selection for systematic practical applications.
Our Kalman SINDy architecture is shown to be competitive with other combinations of data smoothing and system identification techniques, and has a significant advantage in preservation of problem structure and ease of parameter selection.

Model discovery methods are emerging as a critical component in data-driven engineering design and scientific discovery.
Enabled by advancements in computational power, optimization schemes, and machine learning algorithms, such techniques are revolutionizing what can be achieved from sensor measurements deployed in a given system.
Of interest here is the discovery of dynamic models, which can be constructed from a diversity of techniques, including simple regression techniques such as the {\em dynamic mode decomposition} (DMD)~\citep{kutz2016dynamic,ichinaga2024pydmd} to neural networks such as {\em physics-informed neural networks} (PINNs)~\citep{Raissi2019}.
In such models, the objective is to construct a proxy model for the observed measurements which an be used to characterize and reconstruct solutions.
While DMD provides an interpretable model in terms of a modal decomposition, most neural network architectures remain black-box without a clear view of the underlying dynamical processes.
Although the number of techniques available are beyond the scope of this paper to review~\citep{cuomo2022scientific,north2023review}, SINDy is perhaps the leading data-driven model discovery method for interpretable and/or explainable dynamic models as it looks to infer the governing equations underlying the observed data.
As such, it discovers relationships between spatial and/or temporal derivatives, which is the underlying mathematical representation of physics and engineering based systems since the time of Newton.

The SINDy regression architecture seeks to robustly establish relationships between derivatives.
Emerging from~\citep{Brunton2016,brunton2022data}, all variants aim to discover a sparse symbolic representation of an autonomous or controlled system, $\dot x = f(x)$.
A diversity of methodological innovations have been introduced into the SINDy discovery framework to make it robust and stable, including the \emph{weak form} optimization by Messenger and Bortz~\citep{messenger2021bweak,messenger2021weak}.
This approach solves the sparse regression problem after integrating the data over random control volumes, providing a dramatic improvement to the noise robustness of the algorithm.
Weak form optimization may be thought of as a generalization of the integral SINDy~\citep{Schaeffer2017} to PDE-FIND.
Further improvements to noise robustness and limited data may be obtained through ensembling techniques~\citep{Fasel2022}, which use robust statistical bagging to learn inclusion probabilities for the sparse terms $\xi$, similar to Bayesian inference~\citep{gao2022bayesian,gao2023convergence,Hirsh2022}.
Many methodological innovations are integrated in the open-source PySINDy software library~\citep{Kaptanoglu2022}, reducing the barrier to entry when applying these methods to new problems. Additional techniques for learning dynamics from data include PDE-NET~\citep{Long2019,Long2018} and the Bayesian PDE discovery from data~\citep{atkinson2020bayesian}.
Symbolic learning has also been developed, including symbolic learning on graph neural networks~\citep{cranmer2019learning,cranmer2020discovering,sanchez2020learning}.

Kalman smoothing, which this paper integrates with SINDy, has a long history of assimilating measurement data in time series.
From its debut in~\cite{kalman, KalBuc}, engineering practice and design have used it for control and prediction across the real world, e.g. in radar systems, econometric variables, weather prediction, and more.
The family of Kalman methods encompasses both smoothing, after-the-fact techniques, and filtering, real-time updates, that derive from the same assumptions for distributions.
The Kalman smoother can be considered as a best-fit Euler update, the maximum likelihood estimator of Brownian motion, or as the best linear fit of an unknown system.
The best fit/maximum likelihood view extends the classic Kalman updates to a rich family of efficient generalized Kalman smooth algorithms for signals corrupted by outliers, nonlinear models, constraints through side information, and a myriad of other applications, see~\cite{aravkin2017generalized,aravkin2012robust,jonker2019fast}.
In the simplest invocation, the Kalman estimator is determined given only the ratio of measurement noise to the process's underlying stochastic noise.
Fixing both of these parameters allows Kalman methods to also identify the variance of the associated estimator.
Furthermore, a line of research aims to identify parameters purely from data, including \cite{Barratt2020,VanBreugel2020,jonker2020efficient}.
Many methods include their own parameters and are not guaranteed a solution, but are an improvement on the indeterminate nature of direct maximum likelihood or MAP likelihood.

This paper introduces Kalman smoothing as the derivative estimation step in SINDy in distinction with the L1 total variation minimization or Savitzky-Golay smoothers common in application.  It is not the first to combine Kalman methods with SINDy; \cite{rosafalco_ekf-sindy_2024} utilize Ensemble Kalman Filtering (EKF) to identify a partially-known system as a portion of a multi-step method, and \cite{wang_time-variant_2022} apply Kalman filtering to the ODE coefficients as a way of modeling a non-stationary but separable system.
This paper's introduction of Kalman smoothing a continuous process loss for derivative estimation, on the other hand, begins to align the derivative estimation step to the symbolic regression step.
It allows engineering applications to incorporate SINDy estimation with a well-established and familiar data assimilation technique whose noise properties are well understood.

Section two describes the individual methods of SINDy and Kalman smoothing, providing some literature review. In section three, experiments demonstrate the advantages of incorporating Kalman with SINDy. The paper concludes with avenues for future research in section four.

\section{Background}

\subsection{SINDy}
SINDy \citep{Brunton2016} is a family of emerging methods for discovering the underlying dynamics of a system governed by unknown or partially-known \citep{Champion2020} differential equations. It can handle ODEs as well as PDEs~\citep{Rudy2017}, and has been used for protein folding~\citep{Boninsegna2018}, chemical reaction networks~\citep{Hoffmann2019}, plasma physics~\citep{Guan2021}, and more. Most invocations occur through the pysindy Python package, but innovations such as Langevin Regression~\citep{Callaham2021} or \cite{Rudy2019} exist as independent code.

\Figure[t!](topskip=0pt)[width=.4\pdfpagewidth]{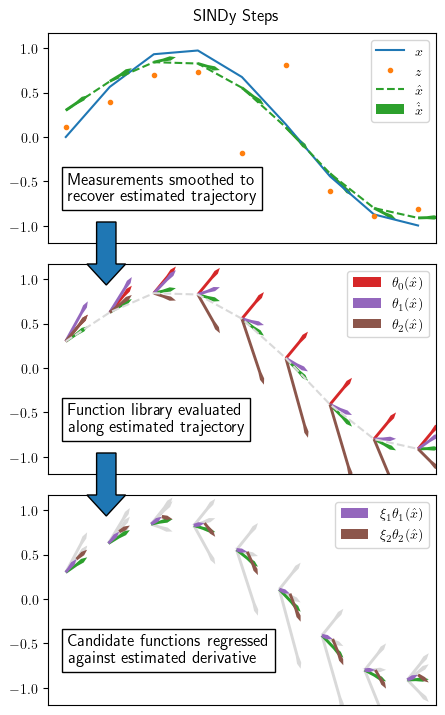}{
    \textbf{The SINDy method, applied to fitting a sinusoid. Taking noisy data, it identifies the model $\dot x = \xi_1\theta_1(x) + \xi_2\theta_2(x)$, rejecting $\theta_0$.}\label{fig:sindy}}

Given some variable of interest $\vec X$ and a library of functions $\mat \Theta$ (including spatial derivatives, when relevant) SINDy seeks to find the coefficients $\mat \Xi$ of the differential equation:

\begin{align}
    \label{eqn:sindy_ode}
    \dot {\mat X} = \Xi {\mat\Theta}({\mat X}),
\end{align}

where
\begin{align*}
    &\mat X \in \R{n \times m}=\vec x(t_1)... \vec x(t_m)\text{: system of $n$ coordinates at $m$ timepoints.}\\
    &\mat \Theta({\mat X}) \in \R{p \times m}\text{: library of $p$ functions evaluated at $m$ timepoints}\\
    &\mat \Xi \in \R{n \times p}\text{: coefficients for $n$ equations of $p$ functions}
\end{align*}
The function library written as a time-independent quantity refers to the collection $\mat \Theta = [\theta_1, \dots \theta_p]^T$, where $\theta_i: \R{n}\rightarrow\R{}$. Examples include the family of all degree-2 polynomials of $n$ inputs, mixed sines and cosines of certain frequencies, or any user-specified family.

The method generally presumes the measurements ($\mat Z$) faithfully reflect system state ($\mat X$) and proceeds in two steps:

\begin{enumerate}
    \item Estimate the time derivatives of the system ${\mat{\widehat{\dot X}}} =\mathrm F(\mat Z)$ for some smoothing function $\mathrm F$.
    \item Choosing a sparse regression method, solve the problem $\underset{\text{sparse } \mat \Xi}{\arg\min} \left\| \mat{\widehat{\dot{X}}} - \mat \Xi \mat \Theta(\mat X) \right\|^2$.
\end{enumerate}
This general process is sketched out in Fig. \ref{fig:sindy}. Researchers have tried a few different methods for calculating the derivatives, broadly grouped into  global methods (e.g. L-1 total variation minimization of \cite{Chartrand2011}) and local methods (e.g. Savitzky-Golay smoothing). Different ways of applying sparsity has attracted more attention, including sequentially thresholding linear regression, nonconvex penalties such as L-0 with a relaxation-based minimization method \citep{zheng2018unified,Champion2020}, an L-0 constraint \citep{Bertsimas2023}, and Bayesian methods for a prior distribution such as spike-slab or regularized horseshoe priors \citep{Hirsh2022,gao2022bayesian}. The latter two papers also demonstrate an interesting line of innovation, eschewing derivatives and using the integral of function library in the loss term. A related approach instead uses the weak form of the differential equation, yielding a solution that {\it is} convex, but which does not provide as straightforward an interpretation of the measurement noise. Most of these methods can benefit from ensembling the data and library terms, as in \cite{Fasel2022}, but others, such as \cite{Kaptanoglu2021} for identifying Galerkin modes of globally stable fluid flows, require a specific form of function library.

This paper seeks to make SINDy more resilient to noise by taking a data assimilation approach. It instead presents the Kalman SINDy steps:
\begin{enumerate}
    \item Estimate the state and time derivatives of the system ${\mat{\widehat{\dot X}}}, \mat{\widehat X} = \mathrm F(\mat X)$ where $\mathrm F$ applies Kalman smoothing.
    \item Choosing a sparse regression method, solve the problem $\underset{\text{sparse } \mat \Xi}{\arg\min} \left\| \mat{\widehat{\dot{X}}} - \mat \Xi \mat \Theta(\mat {\widehat X}) \right\|^2$.
\end{enumerate}

\subsection{Kalman Smoothing}

Kalman filtering and smoothing refers to a group of optimal estimation techniques to assimilate measurement noise to a random process.
Filtering refers to incorporating  new measurements in real-time, while smoothing refers to estimating the underlying state or signal using a complete trajectory of (batch) measurements.
While the processes this paper is concerned with are not random, in the first step of SINDy they are unknown, and so probabilistic language is appropriate.

In adding Kalman smoothing to SINDy, we introduce a distinction between the measurement variables and the state variables of the dynamical system in equation
\ref{eqn:sindy_ode}. As such, the inputs to the problem become $m$ time points of measurements of $k$ variables ($\mat Z\in \R{k\times m}$) and a linear transform from the state to the measurement variables $\mat H \in \R{n \times k}$ describing how the process is measured.

Measurement error is assumed to be normally distributed with $\mat H \mat X - \mat Z \sim \sigma_z \mathcal N(0, \mat R)$ where the covariance matrix $\mat R\in\R{k \times k}$. Measurement regimes where noise is autocorrelated or varies over time can be accomodated by flattening $\mat H \mat X - \mat Z$ and describing $\mat R\in\R{nk \times nk}$.

As a simplifying assumption for experiments in this paper, we use $\mat R=\mat I$. Two parameters are required: $\sigma_z$, the measurement noise standard deviation, and $\sigma_x$, the process velocity standard deviation per unit time. If only point estimates of the state are required, and posterior uncertainty is not, it suffices to use the ratio $\rho = (\sigma_z / \sigma_x)^2$.
\Figure[t!](topskip=0pt)[width=.8\pdfpagewidth]{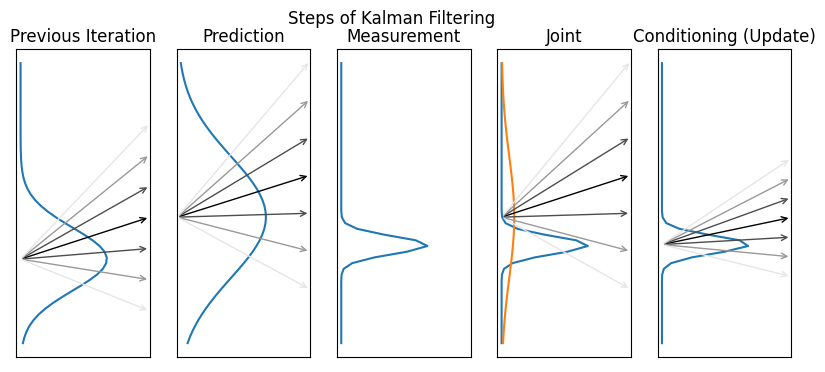}{
    \textbf{Explanatory depiction of Kalman filtering.
    A previous iteration gives a distribution $p(x_{i-1},\dot x_{i-1})$.
    Multiplication by an update matrix produces the predictions $p(x_i, \dot x_i|x_{i-1}, \dot x_{i-1})$.
    Simultaneously, measurements $z_i$ are taken that, with known measurement noise, give $p(z_i|x_i)$.
    Multiplication gives the joint distribution $p(x_i, \dot x_i, z_i| x_{i-1}, \dot x_{i-1})$, from which the conditional distribution $p(x_i, \dot x_i|z_i, x_{i-1}, \dot x_{i-1})$ can be calculated, shown in \cite{Eaton2007}}\label{fig:kalman}}

Each process is assumed to have an independent, Brownian velocity. This leads to Kalman smoothing estimator:
\begin{align}
    \label{eqn:kalman}
    \underset{X, \dot X}{\arg\min}{\|\mat H \mat X - \mat Z\|_{\mat R^{-1}}}^2 + \rho {\|\mat G [\mat {\dot X}, \mat X]\|_{\mat Q^{-1}}}^2.
\end{align}
Here, $\mat G$ is a linear transform to separate $[\mat{\dot X}, \mat X]$ into independent, mean-zero increments, and $\mat Q$ is the covariance of those increments. A graphic displaying Kalman filtering is shown in Fig. \ref{fig:kalman}. To illustrate the ideas, the figure presents step-by-step  filtering updates; however,  batch smoothing is used for the model discovery applications presented in the experiments.

We use the generalized cross validation of \cite{Barratt2020} to choose $\rho$. This strategy chooses $\rho$ in order to minimize the loss on a witheld set of data. While the algorithm described in that paper is not guaranteed to find a minimum, heuristic experience has shown that the longer the trajectory, the more likely their algorithm will succeed. The experiments in the next section show that this strategy works reasonably well.

The generalized cross validation approach of \cite{Barratt2020} witholds some measurement points in order to find the values of $\mat H, \mat R, \mat G,$ and $\mat Q$ that produce estimates $\hat{\mat X}$ that fit witheld data most accurately. This powerful approach can apply to all linear systems but comes with the burdens of nonconvexity. In our work, we presume to know the measurement parameters and most of the process parameters - after all, ``position is the integral of velocity", implies certain constraints on $\mat G$ and $\mat Q$. The method accomodates these constraints via specification of a prox function.

\section{Experiments}

We seek to evaluate Kalman smoothing as a step in SINDy in comparison to other noise-mitigation innovations. We simulate eight dynamical systems\footnote{Cubic Harmonic Oscillator, Duffing, Hopf, Lotka-Volterra, Rossler, Simple Harmonic Oscillator (SHO), Van Der Pol Oscillator, and Lorenz-63} with noisy measurements across a variety of initial conditions, discovering ODEs from SINDy with different smoothing methods. 
The trials are run across a range of durations and relative noise levels, calculated as the noise-to-signal ratio of the measurement variance with the system's mean squared value.
To compare the methods, we then integrate the discovered equations and observe how well they preserve the system's structure as well as directly comparing the coefficients through the F1 score and mean absolute error (MAE).

We compare the results of SINDy with Kalman smoothing and the hyperparameter optimization of \cite{Barratt2020} in comparison with alternative smoothing methods: L-1 total variation minimization and Savitzky-Golay. The latter smoothing methods have been modified to pass the not just the smoothed derivatives $\mat{\widehat{\dot X}}$, but also the smoothed position estimates $\mat {\widehat X}$ to the second step of SINDy. They also each require a parameter: TV requires a coefficient for the L-1 regularizer and Savitzky-Golay requires a smoothing window. These are gridsearched over a wide range, although it is worth noting that choosing the gridsearched optimum requires knowledge of the true system, in distinction to the hyperparameter optimization method used for Kalman smoothing.
We also compare with a gridsearched Kalman smoothing to directly evaluate the efficacy of the generalized cross-validation hyperparameter selection.

Beyond the differentiation step, the SINDy models also specify a function library and optimizer. The feature library used for all experiments was cubic polynomials, including mixed terms and a constant term. The optimizer was the mixed-integer SINDy optimizer of \cite{Bertsimas2023}, configured with the correct number of nonzero terms a priori, and ensembled over 20 models each trained on 60\% of the data. Presenting SINDy with the known number of nonzero coefficients is an attempt to present a best case, where we can ameliorate any interaction between the smoothing method and sparsification parameters. A full list of ODE, simulation, and experimental parameters are shown in the Appendix, tables \ref{tab:ODEs} and \ref{tab:exp-params}.

Methods can be compared in several ways: by the coefficients of the equations they discover, by their accuracy in forecasting derivatives, and how well the discovered system recreates observed dynamics in simulation. As \cite{Gilpin2023} notes, there are many metrics for scoring dynamical system discovery, and the merit of a metric depends upon both the use case and whether the trajectory considered is one of importance. For instance, in controls engineering, the local derivative and very short-term forecasting is the primary imperative. On the other hand, for reduced-order PDE models, recreating larger-scale phenomena in simulation may be more important. Finally, in high-dimensional network dynamics, the accuracy of identifying connectivity, as measured by coefficient F1 score, is most important.

As the coefficient metrics are the most straightforwards, and we compare methods by F1 score and Mean Absolute Error as the duration of training data increases, and separately, as the measurement noise increases. We also visually evaluate how well the discovered ODEs, simulated from random initial conditions in a test set, track the true data and display relevant behavior.

\subsection{Running Experiments}

In a desire to make the experiments not just reproducible, but also reusable, we have separated the method, experiment, and experiment runner into separate packages.
Methodological improvements include adding Kalman smoothing and a entry point for hyperparameter optimization to the \verb|derivative| package, as well as an API for returning not just the derivatives, but the smoothed coordinates themselves (employed for Kalman and Total Variation).
In \verb|pysindy|, we enabled incorporating the smoothed coordinates into successive problem steps.
It should be noted that previous experiments using \verb|pysindy|'s derivative estimation would re-use the noisy coordinates in function library evaluation.

Within \verb|pysindy|, we redefined ensembling in a more flexible way to apply to a greater variety of underlying optimizers, such as the MIOSR one used in these experiments. The standardization of interfaces allows us to compose SINDy experiments in the \verb|pysindy-experiments| package~\citep{pysindy-experiments-2024}. It allows a standard API to specify data generation, model building, and evaluation.

Finally, in order to make it easier to collaborate and reproduce experiments, we expanded the \verb|mitosis| package \citep{Stevens-Haas_mitosis}. This package allows specifying experiment parameters and groups in a declarative manner, which leads to more readable diffs. It also pins reproducibility information for any experiment run. Further reproducibility info is in the Appendix.

\Figure[t!][width=.8\pdfpagewidth]{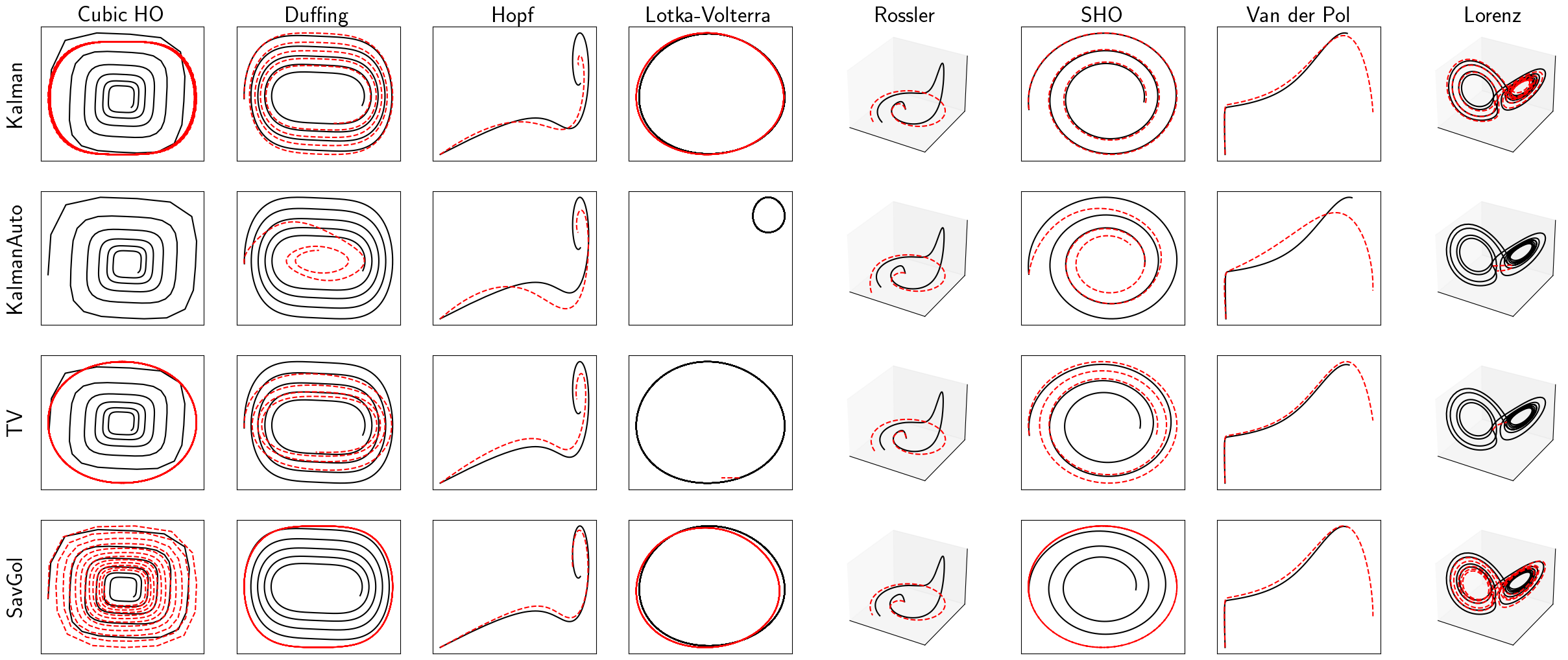}
{ \textbf{The simulation of discovered models compared to test data. Kalman appears better for half of eight ODEs. It represents the essential behavior of more ODEs than TV and Savitzky-Golay. Kalman with auto-hyperparameter selection performs similarly to total variation on a gridsearch. 10\% relative noise, 8 seconds of data.}\label{fig:test}}

\subsection{Results}

We find that SINDy with Kalman smoothing recovers the problem structure in application as well or better than competing methods.
Models discovered in this manner track the essential dynamics in most cases. SINDy with Kalman hyperparameter optimization tends to perform worse than that with Savitzky-Golay, but on par with Total Variation gridsearched optima, and is itself outperformed only slightly by the Kalman gridsearched optima.
While hyperparameter optimization imposes some runtime cost, it does not require access to the true data, making those results all the more inspiring for field use cases.
Simulations of discovered models across all ODEs and methods are shown in Fig. \ref{fig:test}

Surprisingly, methods that smooth better, as shown in Fig. \ref{fig:train} do not necessarily recreate the essential dynamics better in simulation. As a case in point, the Kalman-smoothed training data itself does not seem as accurate as data smoothed by L-1 Total Variation in the the SHO trial.

\Figure[t!][width=.8\pdfpagewidth]{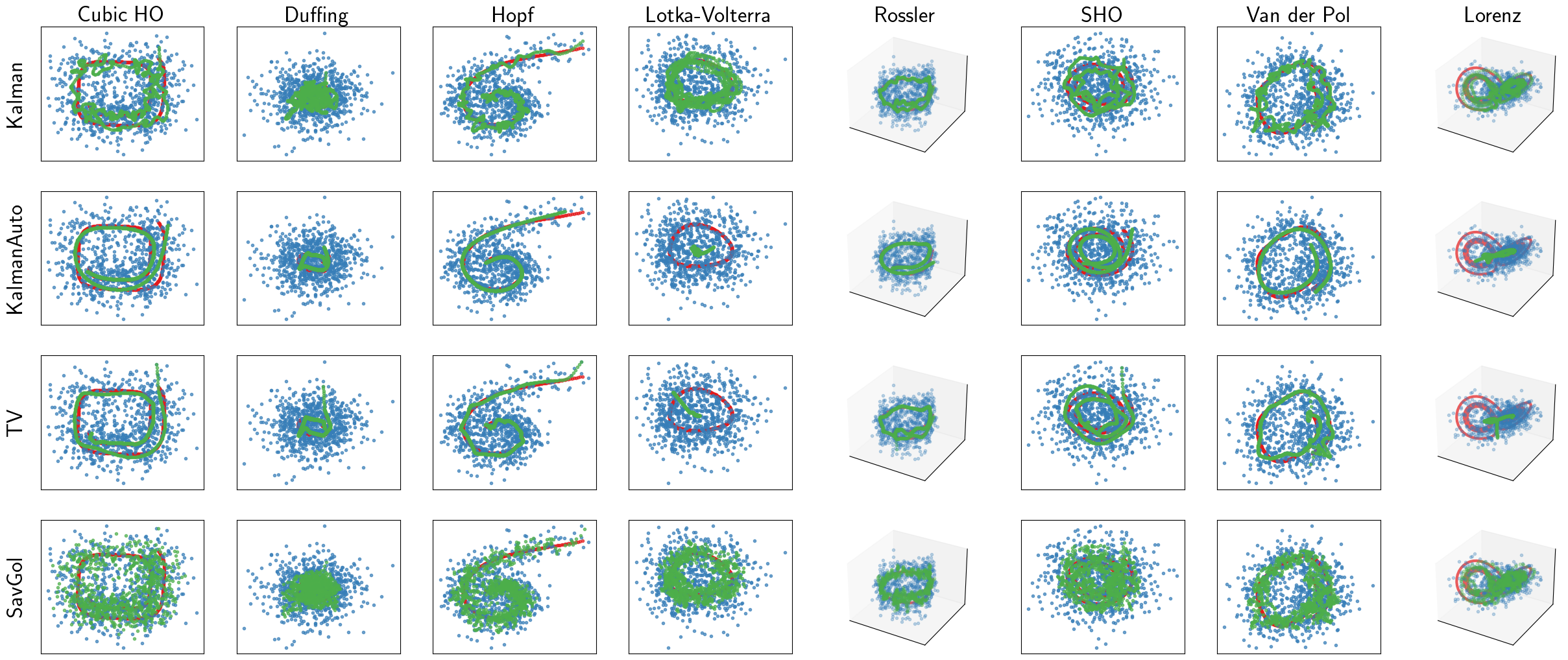}
{ \textbf{The smoothing of training data, performed by different differentiation methods prior to SINDy fit. It does not appear to be the case that a more visually accurate smoothing yields a model that behaves more correctly in simulation. Nevertheless, as Fig. \ref{fig:test} shows, Kalman-smoothed trajectories lead to better models in simulation. 10\% relative noise, 8 seconds of data.}\label{fig:train}}

\begin{figure*}[t!]%
          \centering%
          \vspace*{5pt}%
          \vspace*{0pt}
          \includegraphics[width=\columnwidth]{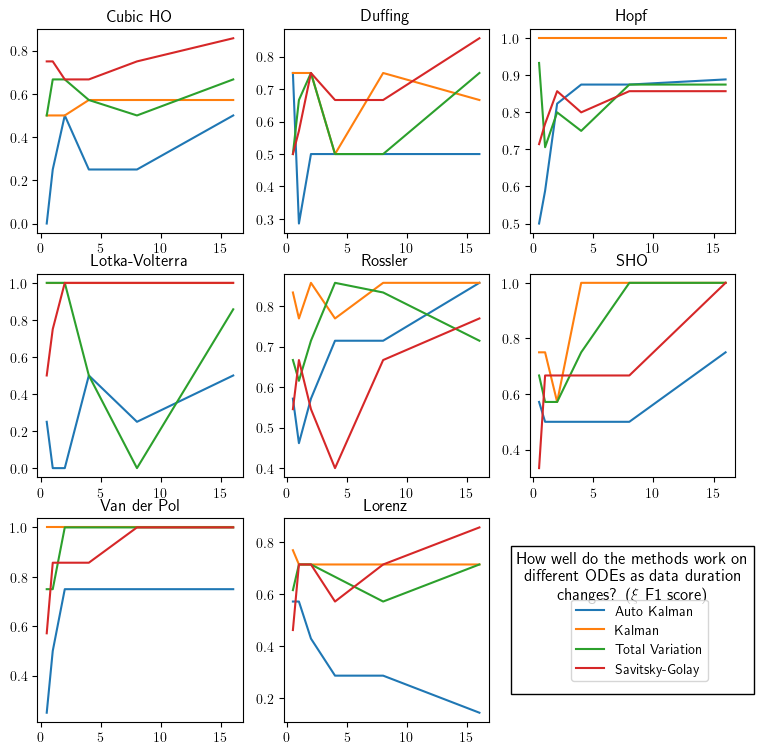}%
          \includegraphics[width=\columnwidth]{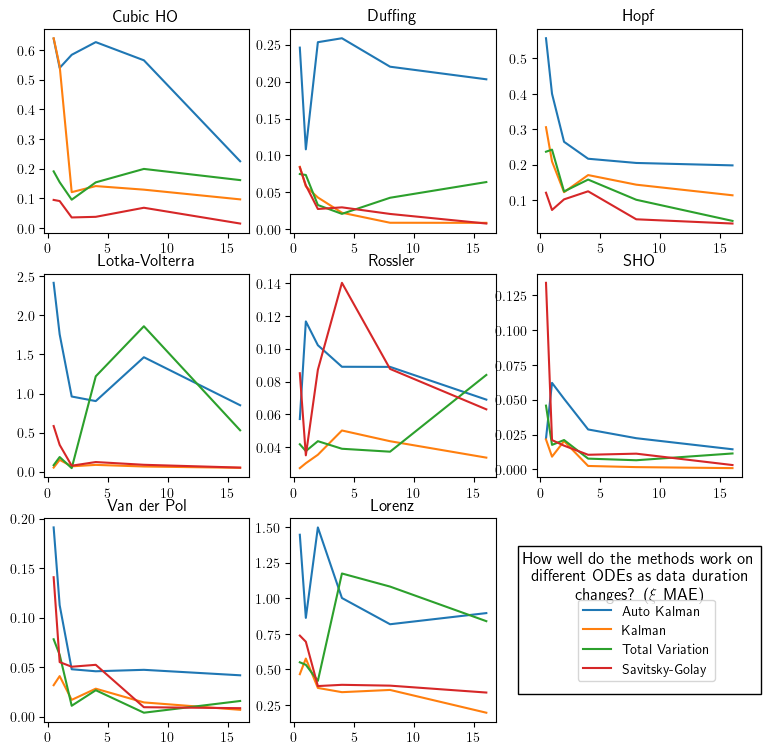}%
          \vspace*{0pt}
          \caption{\textbf{How well different smoothing methods in SINDy recover the ODE coefficients as data duration increases}}%
          \label{fig:time}
          \vspace*{0pt}
\end{figure*}

\begin{figure*}[t!]%
          \centering%
          \vspace*{5pt}%
          \vspace*{0pt}
          \includegraphics[width=\columnwidth]{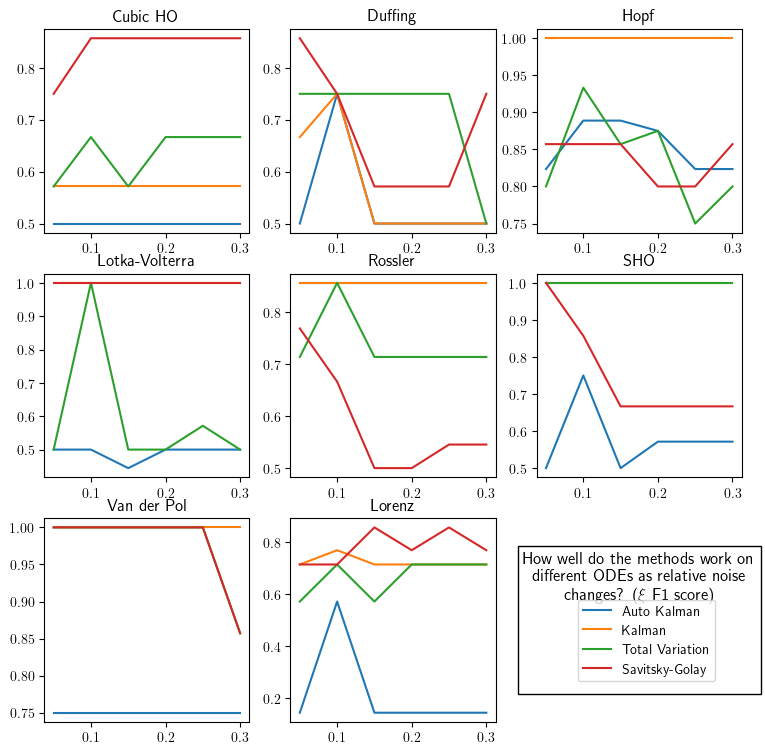}%
          \includegraphics[width=\columnwidth]{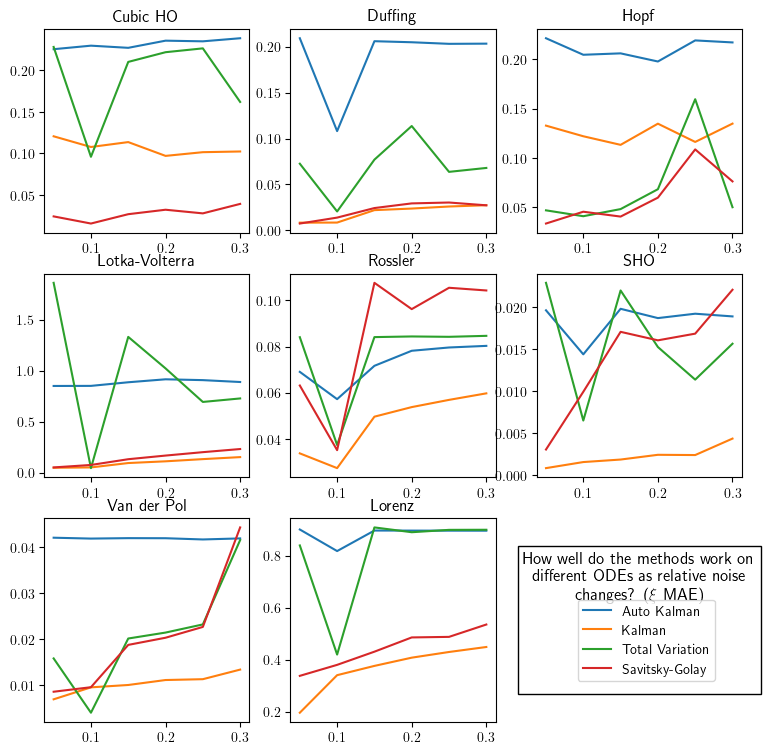}%
          \vspace*{0pt}
          \caption{\textbf{How well different smoothing methods in SINDy recover the ODE coefficients as noise increases}}%
          \label{fig:noise}
          \vspace*{0pt}
\end{figure*}

Even more surprisingly, despite performing well in simulation, SINDy with Kalman hyperparameter optimization performs middlingly in the coefficient metrics. If there's anything that appears consistent about Kalman with GCV, it is that, with a long enough duration, performance appears insensitive regardless of noise, as shown in Fig. \ref{fig:noise}. Across the range of noise levels sampled, either Savitzky-Golay or Kalman (gridsearched) perform the best, depending upon system. As expected, Kalman (gridsearched) always outperforms Kalman GCV, but it is interesting to note that at some durations and noise levels Kalman GCV occasionally outperforms Savitzky-Golay (e.g. Rossler, Hopf).  Coefficient metrics by data duration is shown in Fig. \ref{fig:time}.

Generally, MAE seems to provide a better indication of which method will perform better in simulation than F1 score. Nevertheless, there are cases where the MAE scores of different methods do not indicate which method performs better in simulation, and where effective smoothing does not predict effective system recovery. As one case in point, Kalman GCV and Total Variation smoothing appears most accurate for the Hopf system in Fig. \ref{fig:train}.  However, the coefficient metrics show that either Kalman or Savitzky-Golay recovered the system equations better, and Fig. \ref{fig:test} shows that Savitzky-Golay reconstructed the dynamics more accurately.  Similarly, methods have a wide range of performance on MAE and F1 score on the Rossler system, despite all simulations missing the chaotic behavior.

\section{Conclusion}
This paper has demonstrated that Kalman smoothing is a useful addition to SINDy. It makes the method more generally applicable across domains. The Kalman smoother behaves optimally for the simplest systems and provides a familiar process to the controls engineering community. It also appears to perform better at preserving global system structure in simulation. Incorporating the GCV hyperparameter optimization of \cite{Barratt2020} may not recover the best model, but it allows one to at least recover useful models without relying on an accurate parametrization a priori, particularly if substantial training data exists.  However, ``best model" means different things for different use cases. For uncovering connections between variables, such as in neural activity or chemical reaction networks \cite{Hoffmann2019}, the performance on coefficient F1 metrics indicates that more accurate parameter tuning is essential.

The field is rife with a diversity of follow up studies. Firstly, since Kalman smoothing and SINDy regression loss terms both accommodate variable timesteps, a natural innovation is to combine the two into a single optimization problem. \cite{Hirsh2022} and \cite{Rudy2019} introduce a single-step optimization, but do not evaluate their single step methods in comparison to the mathematically nearest two-step SINDy. As a result, it is difficult to evaluate that aspect of their innovations in isolation.
Following this line of inquiry, producing a single-step SINDy that utilizes Kalman loss could allow a more clear trade-off between measurement noise or the coefficient sparsity.

In parallel, more could be done to give hyper-parameter optimization  access to the terms in the SINDy expression. The method of \cite{Barratt2020} is supremely general, with no intrinsic understanding of the process variance or measurement noise. Applying that knowledge to their prox-gradient method was part of this paper. However, the method is nonconvex, which became problematic with the restriction to the scalar $\rho$ in equation \ref{eqn:kalman}. Moreover, on the path to a single-step SINDy lies an opportunity to use knowledge of the ODE terms in hyper-parameter optimization. In a related note, \cite{VanBreugel2020} also provide hyper-parameter estimation techniques for Savitzky-Golay that could be evaluated in the experiments of this paper.

Kalman SINDy could also be more directly compared to Weak SINDy \citep{messenger2021bweak}, which aims at the same goal of reducing the sensitivity to noise. However, it's implementation in pysindy does not allow for simulation, and modifying existing code to provide that comparison is an investigation in its own right and a necessary next step.

Finally, the interpretation of Kalman smoothing as the maximum likelihood estimator for Brownian motion suggests that it could inform the attempts at a stochastic SINDy, for which \cite{Callaham2021} and \cite{Boninsegna2018} have made the first steps. Stochastic SINDy in those cases aimed for the use cases of noise arising from PDE discretization and inherently statistical mechanics, but the general formulation also has use in any randomly-forced system, such as HVAC controls for a building or mapping a limited part of a chemical reaction network.

\appendix[\break Reproducibility]
\section{}
This paper is built from https://github.com/Jacob-Stevens-Haas/Kalman-SINDy-paper.
To run the experiments, install the package located in images/gen\_image and run the commands in images/gen\_image/run\_exps.sh. Each experiment trial will generate a pseudorandom hex key for reproducibility. To build the final figures, edit images/gen\_image/composite\_plots.py with the keys to the experimental results and run it.

While the exact parametrization is in the experimental configuration and package defaults, it is recreated here in Tables \ref{tab:ODEs} and \ref{tab:exp-params}.

\begin{table*}[t!]
    \centering
    \begin{tabular}{|c|c|c|c|}
        \hline
        System & ODE & Experiment parameters & $x_0$ mean\\
        \hline\hline
        Linear Damped Oscillator
            & $\dot {\vec x} = \left[\begin{matrix}-\alpha & \beta \\ -\beta & -\alpha\end{matrix}\right] \vec x$ & (0,0)
            & $\alpha = .1$, $\beta=2 $\\&&&\\
        Lorenz
            & $\dot {\vec x} = \left[\begin{matrix}
                \sigma (x_2 - x_1)\\
                x_1(\rho - x_3) - x_2\\
                x_1 x_2 - \beta x_3
            \end{matrix}\right]$
            & $\sigma=10$, $\rho=28$, $\beta=8/3$
            & (0, 0, 15)\\&&&\\
        Cubic Damped Oscillator
            & $\dot {\vec x} = \left[\begin{matrix}-\alpha & \beta \\ -\beta & -\alpha\end{matrix}\right] \left[\begin{matrix}x_1^3\\x_2^3\end{matrix}\right]$
            & $\alpha = .1$, $\beta=2 $
            & (0,0) \\&&&\\
        Duffing
            & $\dot {\vec x} = \left[\begin{matrix}x_2 \\ -\alpha x_2 - \beta x_1 -\gamma x_1^3\end{matrix}\right]$
            & $\alpha = .2$, $\beta=.05 $, $\gamma=1$
            & (0,0)\\&&&\\
        Hopf
            & $\dot {\vec x} = \left[\begin{matrix}
                -\alpha x_1 -\beta x_2 - \gamma x_1(x_1 ^2 + x_2^2) \\
                \beta x_1 - \alpha x_2 -\gamma x_2(x_1 ^2 + x_2^2)
            \end{matrix}\right]$
            & $\alpha = .05$, $\beta=1 $, $\gamma=1$
            & (0,0)\\&&&\\
        Lotka-Volterra
            & $\dot {\vec x} = \left[\begin{matrix}
                \alpha x_1 - \beta x_1 x_2 \\
                \beta x_1  x_2 - 2 \alpha x_2
            \end{matrix}\right]$
            & $\alpha = 5$, $\beta=1 $
            & $(5, 5)$ \\&&&\\
        Rossler
            & $\dot {\vec x} = \left[\begin{matrix}
                -x_2 - x_3 \\
                x_1 + \alpha x_2\\
                \beta + (x_1 - \gamma) x_3
            \end{matrix}\right]$
            & $\alpha = .2$, $\beta=.2 $, $\gamma=5.7$
            & (0,0,0) \\&&&\\
        Van der Pol Oscillator
            & $\dot {\vec x} = \left[\begin{matrix}
                x_2 \\
                \alpha(1 - x_1^2) x_2 - x_1
            \end{matrix}\right]$
            & $\alpha = .5$
            & (0,0)\\
        \hline
    \end{tabular}
    \caption{The parametrization of ODEs used in these experiments. Mostly from defaults in the pysindy package.}
    \label{tab:ODEs}
\end{table*}

\begin{table}
    \begin{tabular*}{\columnwidth}{c c}
        Parameter & Value\\
        \hline \\
        \multicolumn{2}{c}{\it Simulated Data} \\
        Number of trajectories & 10\\
        Initial Condition ($x_0$) variance & 9\\
        Initial Condition ($x_0$) distribution & Normal*\\
        Measurement error elative noise (default) & 10\% \\
        Trajectory duration (default) & 16 \\
        Measurement interval & 0.01\\
        Random seed & 19 \\
        \multicolumn{2}{c}{\it SINDy model} \\
        Feature Library ($\Theta$) & Polynomials to degree 3\\
        Optimizer & Mixed Integer Optimizer \\
        L2 regularization (coefficitents) ($\alpha$) & 0.01 \\
        Target sparsity & (true value from equation) \\
        Unbiasing & Yes \\
        Feature normalization & No \\
        Ensembling & data bagging \\
        Number of bags & 20 \\
        \multicolumn{2}{c}{\it Experiment} \\
        Trajectory duration (grid) & 0.5, 1, 2, 4, 8, 16\\
        Relative noise (grid) & 0.05, 0.1, 0.15, 0.2, 0.25, 0.3\\
        Measurement:Process variance (Kalman grid) & 1e-4, 1e-3, 1e-2, 1e-1, 1 \\
        L1 regularization (derivative) (TV grid) & 1e-4, 1e-3, 1e-2, 1e-1, 1 \\
        Window length (Savitzky-Golay grid) & 5, 8, 12, 15\\
        \hline
    \end{tabular*}
    \caption{Parametrization of data, SINDy models, and experiments conducted.\\
    *Lotka-Volterra uses a gamma distribution, rather than normal, in order to enforce nonnegativity.}
    \label{tab:exp-params}
\end{table}

\bibliography{main}

\begin{IEEEbiography}
[{\includegraphics[width=1in,height=1.25in,clip,keepaspectratio]{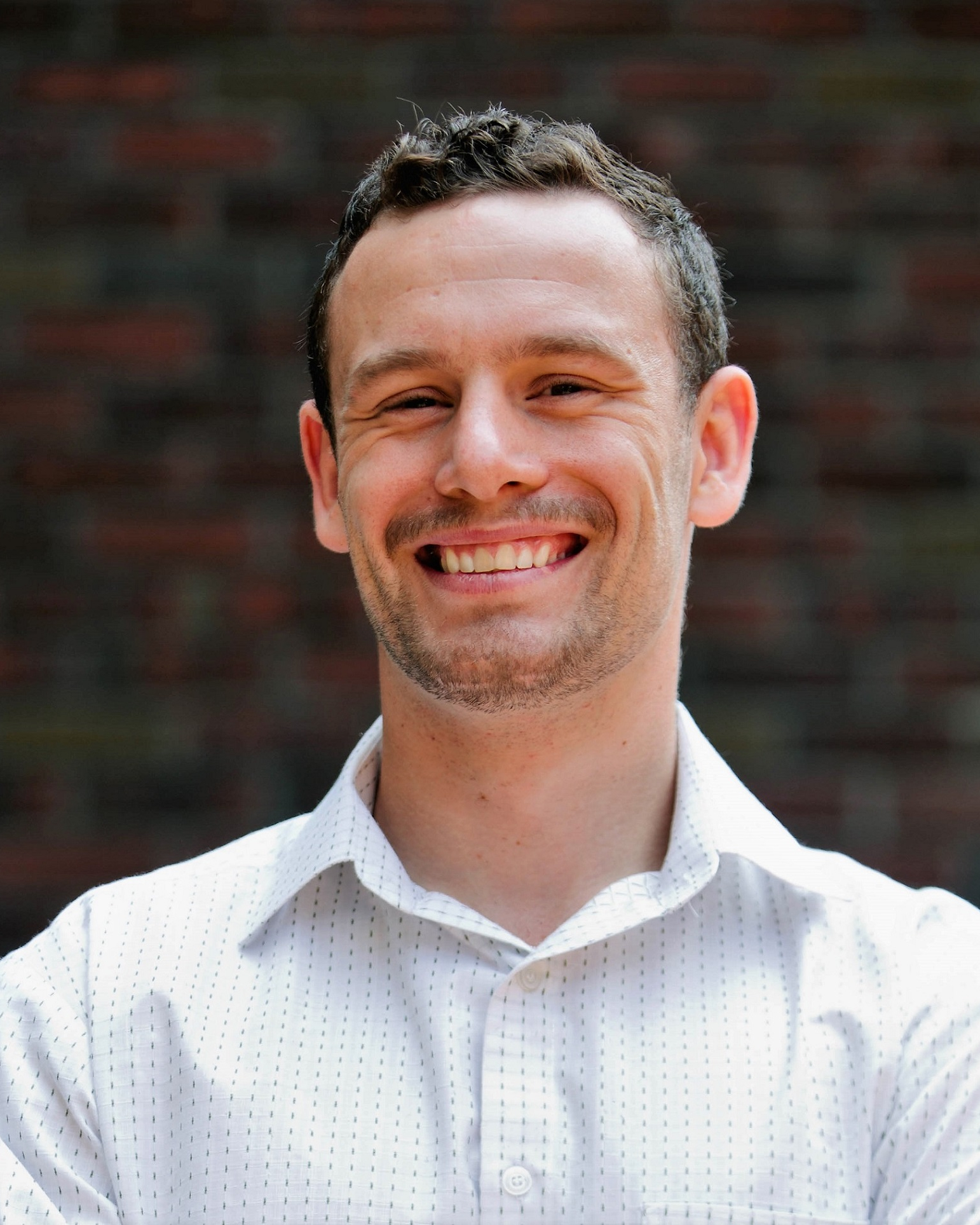}}]
{Jake M. Stevens-Haas} received the B.S. degree in Quantitative Economics from the United States Naval Academy in 2011.

From 2011 to 2017, he was Auxiliaries officer on the USS Momsen (DDG-92) and Weapons and Tactics instructor for Destroyer Squadron NINE. He is currently a PhD. Candidate in Applied Mathematics at University of Washington.

Mr. Stevens-Haas is the lead maintainer of \verb|pysindy|, the python package for Sparse Identification of Nonlinear Dynamics, as well as various patches to various data science Python packages.
\end{IEEEbiography}

\begin{IEEEbiography}
[{\includegraphics[width=1in,height=1.25in,clip,keepaspectratio]{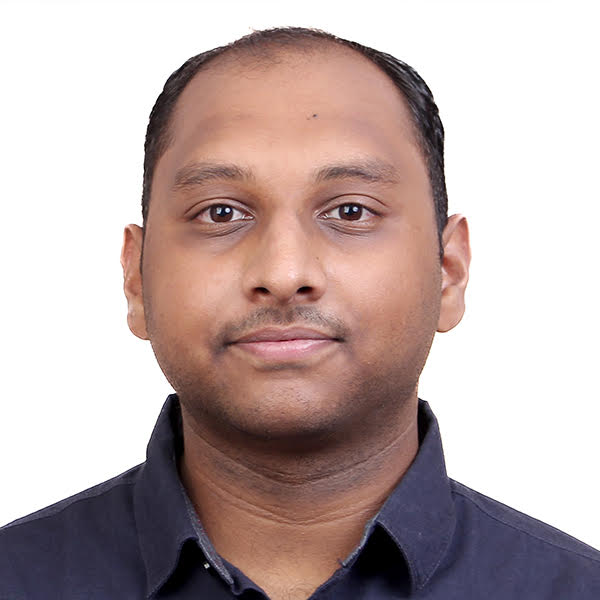}}]
{Yash Bhangale}, born in Jalgaon, Maharashtra, India on May 6th 1999. He received B. Tech. degree in mechanical engineering from the MIT World Peace University at Pune, Maharashtra, India in 2021, is currently pursuing M.S. in mechanical engineering at the University of Washington, Seattle, Washington, United States. His research interests include Dynamical Systems, Machine Learning, System Identification and Differential Equations. 
\end{IEEEbiography}

\begin{IEEEbiography}
[{\includegraphics[width=1in,height=1.25in,clip,keepaspectratio]{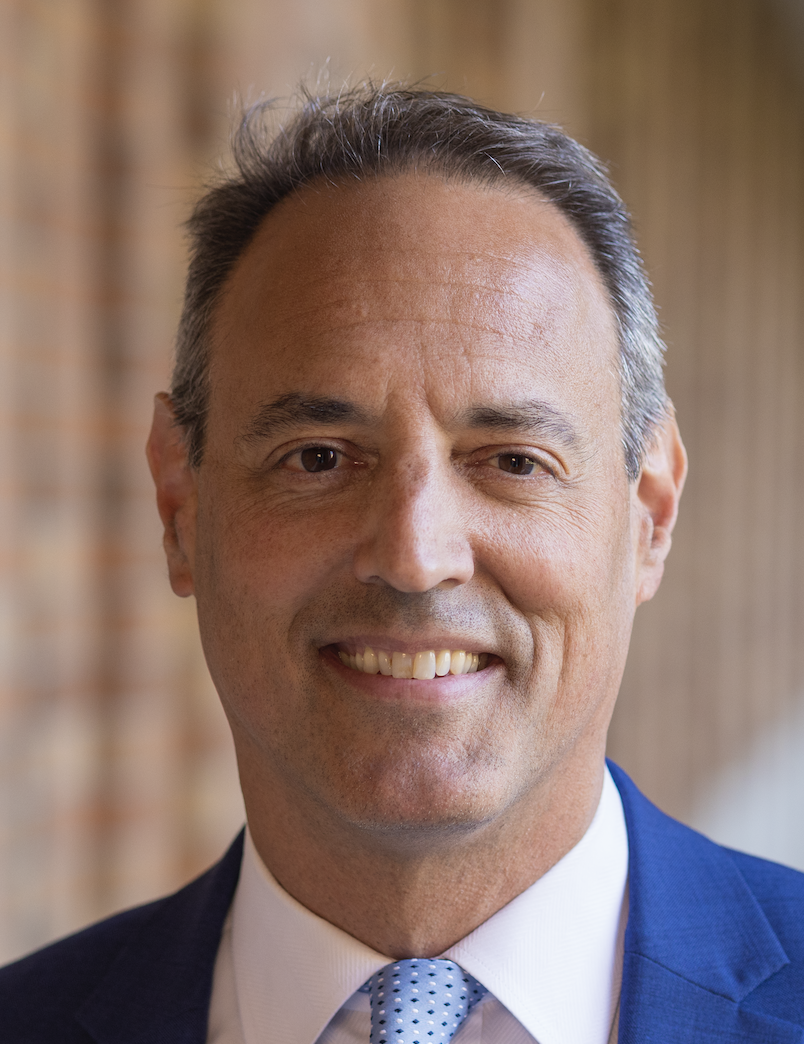}}]
{J. Nathan Kutz} (Senior Member, IEEE) received the
B.S. degree in physics and mathematics from
the University of Washington, Seattle, WA, USA,
in 1990, and the Ph.D. degree in applied mathe-
matics from Northwestern University, Evanston,
IL, USA, in 1994. He is currently Director of the AI Institute in Dynamics Systems, 
a Professor of applied mathematics and electrical engineering, and a Senior Data Science Fellow with the eScience Institute, University of
Washington.
\end{IEEEbiography}

\begin{IEEEbiography}
[{\includegraphics[width=1in,height=1.25in,clip,keepaspectratio]{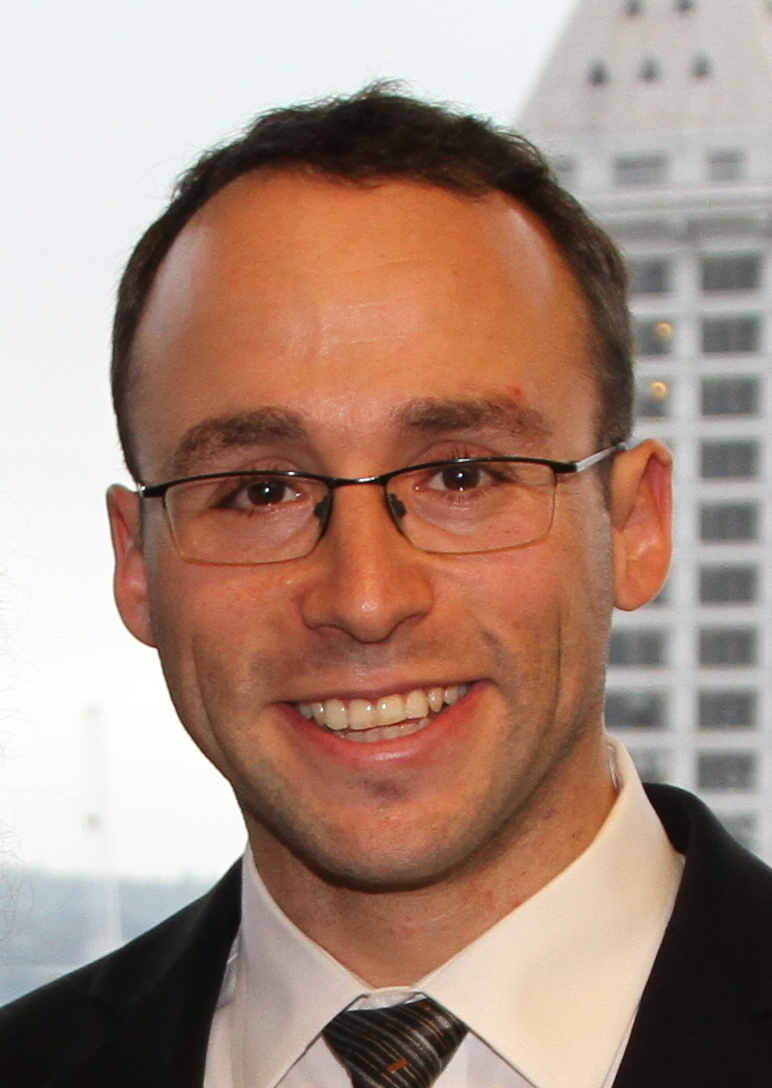}}]
{Aleksandr Aravkin}  received B.S. degrees in Mathematics and Computer Science in 2004, an M.S. in Statistics and a Ph.D. in Mathematics in 2010, from the University of Washington. In 2015, Dr. Aravkin rejoined UW, where he is currently Associate Professor of Applied Mathematics, Adjunct Associate Professor in Statistics, Computer Science, Mathematics, and Health Metrics sciences, Data Science Fellow at the eScience Institute and Director of Mathematical Sciences at the Institute for Health Metrics and Evaluation
\end{IEEEbiography}

\EOD
\end{document}